\documentclass{commat}

%%% AUTHOR'S DEFINITIONS %%%
\theoremstyle{definition}

\DeclareMathOperator{\Alt}{Alt}
\DeclareMathOperator{\Ind}{Ind}
\DeclareMathOperator{\lcm}{lcm}

\title{%
    Structure of finite groups with restrictions on the set of conjugacy classes sizes
    }

\author{%
    Ilya Gorshkov
    }

\affiliation{%
    \address{%
    Sobolev Institute of Mathematics,
        }
    \email{%
    ilygor8@gmail.com
    }}

\abstract{%
    Let $N(G)$ be the set of conjugacy classes sizes of $G$. We prove that if $N(G)=\Omega\times \{1,n\}$ for specific set $\Omega$ of integers, then $G\simeq A\times B$ where $N(A)=\Omega$, $N(B)=\{1,n\}$, and $n$ is a power of prime.}



\VOLUME{32}
\NUMBER{1}
\YEAR{2024}
\firstpage{63}
\DOI{https://doi.org/10.46298/cm.9722}

\begin{paper}

\section*{Introduction}

Let $A, B$ be finite groups and $G := A\times B$. It is easy to check that $N(G)=N(A)\times N(B)$. We are interested in the converse of this assertion.

\begin{question}
Let $G$ be a group such that $N(G)=\Omega\times\Delta$. Which $\Delta$ and $\Omega$ guarantee that $G\simeq A\times B$, where $A$ and $B$ are subgroups such that $N(A)=\Omega$ and $N(B)=\Delta$?
\end{question}

A. Camina proved in \cite{Camin} that, if $N(G)=\{1,p^m\}\times\{1,q^n\}$, where $p$ and $q$ are distinct primes, then $G$ is nilpotent.  In particular, $G=P\times Q$ for a Sylow $p$-subgroup $P$ and a Sylow $q$-subgroup $Q$. Later A. Beltran and M. J. Felipe (see \cite{BF1} and \cite{BF2}) proved a more general result asserting that, if $N(G)=\{1,m\}\times\{1, n\}$, where $m$ and $n$ are positive coprime integers, then $G$ is nilpotent, $n=p^a$ and
$m=q^b$ for some distinct primes $p$ and $q$.

In \cite{SJ}, C. Shao and Q. Jiang showed that if $N(G)=\{1, m_1, m_2\}\times\{1, m_3\}$, where $m_1, m_2, m_3$ are positive integers such that $m_1$ and
$m_2$ do not divide each other and $m_1m_2$ is coprime to $m_3$, then
$G\simeq A\times B$, where $A$ and $B$ are such that $N(A)=\{1,m_1,m_2\}$ and $N(B)=\{1, m_3\}$. In all these cases, the sets of prime divisors of the orders of $A$ and $B$ do not intersect. It was proved in \cite{GorA5A5} that if $N(G)=N(\Alt_5)\times N(\Alt_5)$ and $Z(G)=1$ then $G\simeq \Alt_5\times \Alt_5$.

In \cite{Gorthom1} a directed graph was introduced on the set $N(G)\setminus\{1\}$.
Given $\Theta\subseteq \mathbb{N}$, with $ |\Theta|<\infty$, define the directed graph $\Gamma(\Theta)$, with the vertex set $\Theta$ and edges $\overrightarrow{ab}$ whenever $a$ divides $b$. Set $\Gamma(G)=\Gamma(N(G)\setminus\{1\})$.

In this article, the following theorem is proved.

\begin{theorem}
Let $\Omega$ be a set of integers and $\Gamma(\Omega\setminus\{1\})$ be disconnected, and $n$ be a positive integer such that $\gcd(n,\alpha)=1$ for each $\alpha\in \Omega\setminus\{1\}$. Let $G$ be a finite group such that $N(G)=\Omega\times\{1,n\}$. Then $G\simeq A\times B$, where $N(A)=\Omega$, $N(B)=\{1,n\}$ and $n$ is a prime power.
\end{theorem}

\section{Preliminaries}
We fix the following notation: for an integer $k$, denote by $\pi(k)$ the set of  prime divisors of
$k$. If $\Omega$ is a set of integers, denote $\pi(\Omega)=\bigcup_{\alpha\in\Omega}\pi(\alpha)$. For a prime number $r$, denote by $k_r$ the highest power of $r$ dividing $k$. For integers $m_1,\ldots,m_s$, write
$\gcd(m_1, m_2,\ldots,m_s)$ to denote their greatest common divisor, and write $\lcm(m_1, m_2,\ldots,m_s)$ for their least common multiple.

Let $\Omega$ be a set of integers, and order it by the relation of divisibility. The subset of maximal elements is denoted by $\mu(\Omega)$ and the set of minimal elements is denoted by $\nu(\Omega)$.

\begin{definition}
We say that the set $\Omega$ is separated if, for each $\alpha\in \Omega$, there exists $\beta\in \mu(\Omega)$ such that $\alpha$ does not divide $\beta$.
\end{definition}

Let $G$ be a group and take $a\in G$. We denote by $a^G$ the conjugacy class of $G$ containing $a$. If $N$ is a subgroup of $G$, then $\Ind(N,a)=|N|/|C_N(a)|$. Note that $\Ind(G,a)=|a^G|$. Denote by $|G||_p$ the highest power $p^n$ of $p$ such that $N(G)$ contains multiples of $p^n$ while avoiding multiples of $p^{n+1}$. For $\pi\subseteq\pi(G)$ put $|G||_{\pi}=\prod_{p\in \pi}|G||_p$. For brevity, write $|G||$ to mean $|G||_{\pi (G)}$. Observe that $|G||_p$ divides $|G|_p $ for each $p\in\pi(G)$. In general, $|G||_p $ is less than $|G|_p$.

\begin{definition}
We say that a group $G$ satisfies the condition $R(p)$, or that $G$ is an $R(p)$-group, if there exists an integer $\alpha>0$ such that $a_p\in\{1,p^{\alpha}\}$ for each $a\in N(G)$. In that case, we write $G\in R(p)$.
\end{definition}

The set of $R(p)$-groups can be seen as the disjoints of the two subsets $R(p)^*$ and $R(p)^{**}$:

\begin{enumerate}
\item[$a)$]{ $G\in R(p)^*$ if $G\in R(p)$ and contains a $p$-element $h$ such that $\Ind(G,h)_p>1$;}

\item[$b)$]{ $G\in R(p)^{**}$ if $G\in R(p)$ and $\Ind(G,h)_p=1$ for each $p$-element $h\in G$.}
\end{enumerate}

\begin{lemma}[{\cite[Main theorem]{GorStar}}] \label{star}
If $G\in R(p)^*$, then $G$ has a normal $p$-complement.
\end{lemma}

\begin{lemma}[{\cite[Corollary]{GorStar}}] \label{starcor}
If $G\in R(p)^*$ and $P\in Syl_p(G)$, then $Z(P)\leq Z(G)$.
\end{lemma}

\begin{lemma}[{\cite[Lemma 1.4]{GorA2}}] \label{factorKh}
For a finite group $G$, take $K\unlhd G$ and put $\overline{G}= G/K$. Take $x\in G$ and $\overline{x}=xK\in G/K$.
The following claims hold:
\begin{enumerate}
\item[(i)] $|x^K|$ and $|\overline{x}^{\overline{G}}|$ divide $|x^G|$.

\item[(ii)] For neighboring members $L$ and $M$ of a composition series of $G$, with $L<M$, take $x\in M$  and
the image $\widetilde{x}=xL$ of $x$. Then $|\widetilde{x}^S|$ divides $|x^G|$, where $S=M/L$.

\item[(iii)] If $y\in G$ with $xy=yx$ and $(|x|,|y|)=1$, then $C_G(xy)=C_G(x)\cap C_G(y)$.

\item[(iv)] If $(|x|, |K|) = 1$, then $C_{\overline{G}}(\overline{x}) = C_G(x)K/K$.

\item[(v)] $\overline{C_G(x)}\leq C_{\overline{G}}(\overline{x})$.
\end{enumerate}
\end{lemma}

\begin{lemma}[{\cite[Lemma 4]{GorMal}}] \label{hz2}
Let $g\in G$. If each conjugacy class of $G$ contains an element $h$ such that $g\in C_G(h)$ then $g\in Z(G)$.
\end{lemma}

\begin{lemma}[{\cite[Theorem A]{BFMM}}] \label{vse}
Let $G$ be a finite group, and let $p$ and $q$ be distinct primes. Then
some Sylow $p$-subgroup of $G$ commutes with some Sylow $q$-subgroup of $G$ if and only
if the class sizes of the $q$-elements of $G$ are not divisible by $p$ and the class sizes of
the $p$-elements of $G$ are not divisible by $q$.
\end{lemma}

We call a $p$-element $x$ of $G$ $p$-central
if $x\in Z(P)$ for some Sylow $p$-subgroup $P$ of $G$.

\begin{lemma}[{\cite[Theorem B]{NST}}] \label{nst}
Let $G$ be a finite group and $p$ a prime. Suppose that every $p$-element of $G$
is $p$-central. Then
$$O^{p'}(G/O_{p'}(G))= S_1\times \cdots \times S_r \times H,$$
where $H$ has an abelian Sylow $p$-subgroup, $r\geq0$, and $S_i$ is a non-abelian simple group
with either
\begin{enumerate}
\item[(i)] $p = 3$ and: $S_i\simeq Ru$, or $J_4$, or $S_i\simeq {}^2F_4(q_i)', 9 \not|(q_i + 1)$; or

\item[(ii)] $p = 5$ and $S_i\simeq Th$ for all $i$.
\end{enumerate}
\end{lemma}

\begin{lemma}\label{aberp}
If $G\in R^{**}(p)$, then the Sylow $p$-subgroups of $G$ are abelian.
\end{lemma}
\begin{proof}
Note that $R^{**}(p)$-groups satisfy the condition of Lemma \ref{nst}. Hence, if a Sylow $p$-subgroup is non-abelian, then $p\in\{3,5\}$ and $O^{p'}(G/O_{p'}(G))= S_1\times\cdots\times S_r \times H$, where $S_i$ is isomorphic to one of the groups $Ru, J_4, {}^2F_4(q_i)', Th$. Note that if $r>1$, then the group $G$ is not an $R^{**}(p)$ group. It follows from the description of conjugacy class sizes in \cite{AtlasV3} and \cite{2F4} that $S$ contains a $p'$-element $g_1$ such that $1<\Ind(S,g_1)_p<|S| _p$ and a $p'$-element $g_2$ such that $\Ind(S,g_2)_p=|S|_p$. Since $p$ and $|O_{p'}(G)|$ are relatively prime, there exists $g_1'\in G$ such that $g_1'O_{p'}(G)=g_1$ and $\Ind(G,g_1')_p=\Ind(G/O_{p'}(G),g_1)_p$. Let $g_2'\in G$ be such that $g_2'O_{p'}(G)=g_2$. We have $C_G(g_2')O_{p'}/O_{p'}\leq C_{G/O_{p'}(G)}(g_2)$. In particular $\Ind(G,g_2')_p\geq \Ind(G/O_{p'}(G),g_2)_p>\Ind(G,g_1')_p$, contradicting the definition of $R^{**}(p)$-groups.
\end{proof}

\begin{lemma}\label{Rpssfak}
Any $R^{**}(p)$-group  contains at most one non-abelian composition factor whose order is divisible by $p$.
\end{lemma}
\begin{proof}
Let $G$ be an $R^{**}(p)$-group. Lemma \ref{aberp} implies that the Sylow $p$-subgroup of $G$ is abelian. Let $1<G_1<\cdots<G_k=G$ be the chief series. Assume that $G_i/G_{i-1}=H$ is a non-solvable group and the order of $H$ is divisible by $p$. Lemma \ref{factorKh} implies that the conjugacy class sizes of the group $H$ divide the corresponding conjugacy class sizes of $G$. We have $H=S_1\times S_2\times\cdots \times S_t$, where the $S_i$ are isomorphic non-abelian finite simple groups, for $1\leq i\leq t$.

Assume that $|G_{i-1}|$ is divisible by $p$. Let $P\leq G_{i-1}$ be a Sylow $p$-subgroup of $G_{i-1}$. From Frattini's argument, it follows that $N_{G_i}(P)/N_{G_{i-1}}(P)\simeq G_i/G_{i-1}$. Let $\widehat{H}\leq N_{G_i}(P)$ be a subgroup generated by all Sylow $p$-subgroups of $N_{G_i}(P)$.
Since any Sylow $p$-subgroup of $G$ is abelian and $H$ is generated by $p$-elements, we infer that $\widehat{H}G_{i-1}/G_{i-1}=H$ and $\widehat{H}$ centralizes some Sylow $p$-subgroup of the group $G_{i-1}$.

Assume that $g\in G/G_{i-1}$ is a $p$-element acting on $H$ as an outer automorphism. The fact that the Sylow $p$-subgroups of $G$ are abelian implies that $S_j^g=S_j$ for any $1\leq j\leq t$. Assume that $g$ acts non-trivially on $S_j$. Since the Sylow $2$-subgroup of a simple alternating group of degree greater than $5$ is non-abelian and the outer automorphism group of an alternating group is a $2$-group, we obtain that $S_j$ cannot be isomorphic to any of the alternating groups. It follows from \cite{AtlasV3} and Lemma \ref{nst} that $S_j$ cannot be isomorphic to any of the sporadic groups, and therefore $S_j$ is a group of Lie type. In \cite[Theorem 1]{Shen} and in \cite{Walter} it is described when a Sylow $p$-subgroup of a simple group of Lie type is abelian. We can show that $g$ acts on $S_j$ as a field automorphism. It follows from the description of the centralizers of field automorphisms (see \cite[Theorem 4.9.1]{GLS}) that the Sylow $p$-subgroup of $S_j.\langle g \rangle$ is non-abelian, and hence the Sylow $p$-subgroup of $G$ is non-abelian, which is a contradiction. Therefore, it can be considered that $H$ contains a Sylow $p$-subgroup of $G/G_{i-1}$.

Assume that $t>1$. For each $j \in \{1, \ldots, t\}$, there is an element $h_j\in S_j$ such that $\Ind(S_j,h_j)_p=|S_j|_p$. Let $g=h_1\cdots h_t$ and $\widehat{g}\in\widehat{H}$ be some pre-image of the element $g$. Since $\widehat{H}$ centralizes a Sylow $p$-subgroup of $G_{i-1}$ and $\Ind(H,g)$ divides $\Ind(G,\widehat{g})$, we infer that $ \Ind(G,\widehat{g})_p=(\Ind(H,g))_p=|H|_p$. If $t>1$, then $\widehat{H}$ contains an element $\widehat{h_1}$, which is the pre-image of the element $h_1$ such that $1<\Ind(G,\widehat{h_1})_p <|H|_p$. This contradicts the definition of an $R(p)$-group.
\end{proof}

\begin{lemma}[{\cite[Theorem 5.2.3]{Gore}}] \label{Gore5}
Let $A$ be a $\pi(G)'$-group of automorphisms of
an abelian group $G$. Then $G=C_G(A)\times[G,A]$.
\end{lemma}

\begin{lemma}\label{gooddirect}
Let $P\lhd G$ be a Sylow $p$-subgroup of $G$. If $P=A\times B$ with $A$, $B$ normal subgroups of $G$, then $C_G(ab)=C_G(a)\cap C_G(b)$ for any $a\in A $ and $b\in B$.
\end{lemma}
\begin{proof}
The assertion of the lemma follows from the fact that any $p$-element $x$ is uniquely represented as $x=x_ax_b$ where $x_a\in A$ and $x_b\in B$.
\end{proof}

\section{Proof of the Main Theorem}

Let $G$ be as in the hypothesis of the theorem.
We divide the proof of the theorem into 3 propositions. In the preliminary lemma and in Propositions \ref{um} and \ref{tres}, we only use the separation property of the set $\Omega$. The disconnection of the graph $\Gamma(\Omega \setminus\{1\})$ is used only in the proof of Proposition 2.

Note that $G$ has the property $R(p)$ for any $p\in\pi(n)$.
In Propositions \ref{um} and \ref{dois} we prove that $G\not\in R^{**}(p)$. In Proposition \ref{tres} we analyze the case $G\in R^*(p)$ and thus complete the proof of the Main Theorem.

Assume that $G\in R^{**}(p)$ for any $p\in \pi(n)$. In this case, Lemma \ref{aberp} implies that a Sylow $p$-subgroup of $G$ is abelian. It follows from Lemma \ref{vse} that a Hall $\pi(n)$-subgroup exists and is abelian. It follows from the well-known Wielandt theorem that all Hall $\pi(n)$-subgroups are conjugate.

\begin{lemma}\label{kpfak}
The order of any non-abelian composition factor of $G$ is not divisible by $p$.
\end{lemma}
\begin{proof}
Lemma \ref{Rpssfak} implies that $G$ contains at most one non-abelian composition factor $S$ whose order is divisible by $p$. Let $R\lhd G$ be such that $S\leq G/R$. Let $g\in G$ be a $p$-element such that its image $gR\in S$ is not trivial. Let $x\in G$ be an element of minimal order such that $\Ind(G,x)=n$. Since $n$ is minimal with respect to divisibility in $N(G)$, we infer that $|x|=r^{\alpha}$ is a power of a prime $r$. We have that $x$ centralizes Sylow $t$-subgroups for any $t\in\pi(\Omega)$ and, in particular, $x$ centralizes Sylow $t$-subgroups for any $t\in\pi(\Ind(G,g))$. Put $C=C_G(x)$.
Since $S$ is the unique non-abelian composition factor whose order is divisible by $p$, we infer that $S$ is a normal subgroup of $G/R$. Note that $CR/R$ contains Sylow $t$-subgroups of $G/R$ for any $t\in \pi(\Ind(S,\overline{g}))$. Let $T$ be a Sylow $t$-subgroup of $G/R$ for some prime $t\in\pi(G/R)$. Since $S$ is a normal subgroup of $G/R$ we infer that $T\cap S$ is a Sylow $t$-subgroup of $S$. From the fact that finite simple groups do not have Hall $p'$-subgroups for each prime divisor $p$ of its order, we get that group $S$ is generated by its Sylow $t$-subgroups, where $t\in\pi(\Ind(S,\overline{g}))$. Hence $S\leq CR/R$. In particular, $C$ contains a pre-image of the group $S$. Therefore, $C$ contains an $r'$-element $y$ such that $\Ind(C,y)_p>1$. Thus, $\Ind(G,xy)_p>\Ind(G,x)_p$, which is a contradiction.
\end{proof}

Let $O=O_{\pi(n)'}(G)$. Lemma \ref{kpfak} implies that $G/O$ contains a normal $p$-subgroup $\overline{P}$, for some $p\in\pi(n)$. Let $T=O_{\pi(n)}(G/O)$. Assume that $T$ is not a Hall $\pi(n)$-subgroup of $G/O$.  Since a Hall $\pi(n)$-subgroup of $G$ is abelian, we have $T$ is abelian.  The centralizer of $R$ in $G/O$ is a normal subgroup of $G/O$ for each Sylow subgroup $R$ of $T$. For any $g\in G/O$ it follows from the inequality $\Ind(G/O,g)_{p}>1$ that $\Ind(G/O,g)_{\pi(n)}=n$. Using these facts it is easy to obtain a contradiction. Therefore, $G/O$ contains a normal Hall $\pi(n)$-subgroup $\overline{H}$. In particular, we can assume that $\overline{P}$ is a Sylow $p$-subgroup of $G/O$. Let $x\in G$ be an element of minimal order such that $\Ind(G,x)=n$. Since $n$ is minimal by divisibility number of $N(G)$, we infer that $x$ is an element of order $t^{\alpha}$, where $t$ is some prime and $t\not\in \pi(n)$.

\begin{proposition}\label{um}
The image $\overline{x}\in G/O$ of $x$ is trivial.
\end{proposition}
\begin{proof}
Assume that $\overline{x}$ is not trivial.
Lemma \ref{Gore5} implies that $\overline{P}=[\overline{x},\overline{P}]\times C_{\overline{P}}(\overline{x})$. Let $\widetilde{x}\in G/OH$ denote the image of $x$. Since $\pi(G/OH)$ does not contain numbers from the set $\pi(n)=\pi(\Ind(G,x))$, and $\Ind(G/OH,\widetilde{x})$ divides $ \Ind(G,x)$, we infer that $\Ind(G/OH,\widetilde{x})$ is equal to $1$. Hence $\widetilde{x}\in Z(G/OH)$. Thus the subgroup $C_{\overline{P}}(\overline{x})$ is a normal subgroup of $G/O$. Since $p\not\in \pi(G/OH)$, it follows from Maschke's theorem that $C_{\overline{P}}(\overline{x})$ has compliment in $\overline{P}$. In particular $[\overline{x},\overline{P}]$ is a normal subgroup of $G/O$.

Let $P$ be a Sylow $p$-subgroup of $G$, and let $P_1, P_2 \leq P$ be such that $P_1.O/O=[\overline{x},\overline{P}]$ and $P_2O/O=C_{\overline{P}}(\overline{x})$. Since $\lcm(\Ind(O,x), O)=1$, we have $x\in C_G(O)$. The group $C_G(O)$ is a normal subgroup of $G$. We have $C_G(O)O/O\unlhd \overline{G}$ and $\overline{x}\in C_G(O)O/O$. From the fact that $\overline{x}$ acts without fixed points on $[\overline{x},\overline{P}]$  and $[\overline{x},\overline{P}]\unlhd \overline{G}$ it follows that $[\overline{x},\overline{P}]$ is the minimal normal subgroup of $\overline{G}$ which includes $\overline{x}$. In particular $P_1<C_G(O)$.

The fact that the number $\Ind(G,x)_p$ is maximal implies that centralizer of any $t'$-element of $C_G(x)$ contains some Sylow $p$-subgroup of the group $C_G(x)$. Since $O.P_2\unlhd G$, we infer that the centralizer of any $t'$-element from $O$ contains a subgroup conjugate to $P_2$ in $O.P_2$. Suppose there is a $t$-element $y\in O$ such that $\Ind(O.P_2,y)_p>1$. Since $\Ind(G,x)_t=1$, we infer that $C_G(x)$ contains some Sylow $t$-subgroup of $G$. In particular, one can assume that $y\in C_G(x)$. Consider $C_G(xy)$. Let $R$ be a Sylow $p$-subgroup of $G$ such that $\widetilde{R}=R\cap C_G(xy)$ is a Sylow $p$-subgroup of $C_G(xy)$. Since $P_1\leq C_G(O)$, we have $P_1\leq R$ and $P_1\cap C_G(xy)=1$. It follows from the fact that $\Ind(G,x)_p=\Ind(G,xy)_p=|P_1|$ and the fact that $R$ is an abelian group that $R=P_1\times \widetilde{R}$. Note that $\widetilde{R}<C_G(x)$, and hence $\widetilde{R}$ is conjugate to $P_2$ in $C_G(x)$. In particular, $\widetilde{R}$ is conjugate to $P_2$ in $O.P_2$. Therefore $y$ centralizes $\widetilde{R}$ and $\Ind(O.P_2,y)_p=1$, which is a contradiction. Thus any element of $O$ centralizes some Sylow $p$-subgroup. Lemma \ref{hz2} implies that $P_2<C_G(O)$. Thus $G$ contains a normal abelian Hall $\pi(n)$-subgroup $N$.

We have that $P_2$ is a Sylow $p$-subgroup of $C_G(x)$ and $P_2\unlhd C_G(x)$. From the fact that $\Ind(G,x)_p$ is maximal it follows that any $t'$-element centralizes $P_2$. Therefore, we have that $\pi(\Ind(C_G(x),h))\subseteq\{t\}$ for any $h\in P_2$.
Since $C_G(x)$ contains Sylow $r$-subgroups of $G$ for any $r\in\pi(\Omega)$ and a Hall $\pi(n)$-subgroup of $G$ is abelian, we infer that $\pi(\Ind (G,h))\subseteq\{t\}$ for any $h\in P_2$. Let $g\in C_G(x)$ be some $t'$-element. Then $g$ acts on $P_1$, and
\[
\Ind(G,g)_p=\Ind(P_1,g)_p.
\]
Since $\Ind(P_1,g)_p\in\{1,|P_1|\}$, we see that $g$ acts on $P_1$ either trivially or without fixed points.
Note that $x^G=x^N$. Thus, $\Ind(G,a)_{t'}=\Ind(G,b)_{t'}$ for any $a,b\in P_1$ and $\pi(\Ind(G ,c))\subseteq\{t\}$ for any $c\in P_2$. It follows from Lemma \ref{gooddirect} that, for any $p$-element $a$, there exists $k$ such that $\Ind(G,a)_{t'}\in\{1,k\}$. Thus, $\Omega$ contains a number $\alpha$ dividing the index of any $p$-element. Let $h_1\in P_1$ be such that $\Ind(G,h_1)$ is minimal among $\{\Ind(G,g)|g\in P_1\}$, and let $h_2\in P_2$ be such that $\Ind( G,h_2)$ is minimal among $\{\Ind(G,g)|g\in P_2\}$.

Assume that $\Ind(G,h_2)_t\leq \Ind(G,h_1)_t$. Then $\Ind(G, h_2)$ divides $\Ind(G,g)$ for any $p$-element $g$.
Since $\Omega$ is separated, we obtain that $\mu(\Omega)$ contains an element $\beta$ that is not divisible by $\Ind(G,h_2)$. Let $l\in G$ be such that $\Ind(G,l)=\beta$. Since $\Ind(G,h_2)$ does not divide $\beta$, we infer that $C_G(l)$ does not contain $p$-elements. But $\beta$ is not divisible by $p$ and hence $C_G(h)$ contains some Sylow $p$-subgroup, therefore we have a contradiction.

Thus, $\Ind(G,h_2)_t> \Ind(G,h_1)_t$. Since $\Omega$ is separated, we infer that $\mu(\Omega)$ contains an element $\beta$ that is not divisible by $\Ind(G,h_2)$. Let $l\in G$ be such that $\Ind(G,l)=\beta$. Since $\Ind(G,h_2)$ does not divide $\beta$, we see that $|l|$ is divisible by $p$. Further, we have $l=ab$, where $a$ is a $p$-element and $b$ is a $p'$-element. We have that $\Ind(G,a)$ divides $\beta$. From Lemma \ref{gooddirect} and the fact that $\beta$ is not divisible by numbers in $\{\Ind(G,g)|g\in P_2\}$, it follows that $a\in P_1$. It follows from Lemma \ref{gooddirect} that $\Ind(G,abh_2)$ is divisible by $\beta$ and $\Ind(G,h_2)$ contradicting the fact that $\beta$ is maximal in $\Omega$.
\end{proof}

\begin{proposition}\label{dois}
The element $x\not\in O$.
\end{proposition}
\begin{proof}
Assume that $x\in O$. Since $\Ind(G,x)$ is relatively prime to $|O|$, we have $O\leq C_G(x)$. Let $X=\langle x^G\rangle$. The fact that $O$ is a normal subgroup of $G$ implies that $O\leq C_G(X)$. Hence, $X$ is an abelian $t$-subgroup of the group $O$. Let $P$ be a Sylow $p$-subgroup of $G$ such that $P_1=P\cap C_G(x)$ is a Sylow $p$-subgroup of $C_G(x)$. The fact that $x^G=x^{O.H}$ implies that $P_1< C_G(X)$. Thus, any $t'$-element of $O$ centralizes some subgroup conjugate to $P_1$.

Consider $X$ as a $\widetilde{P}=P/P_1$-module. It follows from Lemma \ref{Gore5} that the group $X$ can be represented as $[X,\widetilde{P}]\times C_X(\widetilde{P})$. Since $\widetilde{P}$ acts non-trivially on $X$, we see that $[X,\widetilde{P}]$ is non-trivial. Since for any element $y\in [X,\widetilde{P}]$ we have $\widetilde{P}\cap C_{O.\widetilde{P}}(y)=1$, we infer that $\widetilde{ P}$ acts without fixed points on $[X,\widetilde{P}]$. Hence $\widetilde{P}$ is a cyclic group. %Let $z\in P$ be such that $zP_1$ generates the group $\widetilde{P}$.

Assume that $P_1$ contains an element $f$ such that $\Ind(G,f)>1$.
We will use the fact that graph $\Gamma(\Omega\setminus\{1\})$ is disconnected. Let $\Gamma_1$ be a connected component of the graph $\Gamma(\Omega\setminus\{1\})$ such that $\Ind(G,f)\in \Gamma_1$. Since any $t'$-element centralizes some element from $f^G$, we infer that $\Ind(G,g)_{\pi(n)'}\in \Gamma_1\cup\{1\}$ for any $\{p, t\}'$-element $g$.

Denote by $\Gamma_2$ some connected component of the graph $\Gamma(\Omega\setminus\{1\})$ different from $\Gamma_1$. Let $y\in G$ be such that $\Ind(G,y)\in \Gamma_2$. We have that $y$ is a $\{p,t\}$-element. Assume that $y$ is a $t$-element. Since $\Ind(G,y)_p=0$, we infer that $C_G(y)$ contains a subgroup conjugate to $P_1$, and hence $\Ind(G,y)\in\Gamma_1\cup\{1\}$, deriving a contradiction.

Therefore, if $\Ind(G,g)_p=1$, then $g$ is the product of a $p$-element and an element from the center of $G$. In particular, if $\Ind(G,g)\in \Gamma_2$, where $g$ is an element of primary order, then $\pi(g)=\{p\}$. It also follows from here that $\pi(n)=\{p\}$.

Since $\Gamma_2$ is an arbitrary connected component of $\Gamma(\Omega\setminus\{1\})$ different from $\Gamma_1$, then we can assume that there exists $z\in P$ such that $\Ind(G,z)\in \Gamma_2$. Then $\Ind(G,y)\in\Gamma_2\cup \{1\}$ for any $y\in\langle z\rangle$. This means that $\langle z\rangle\cap P_1\leq Z(G)$. Let $g\in P$ be such that $z\in\langle g\rangle$. Since $\Ind(G,z)$ divides $\Ind(G,g)$ it follows that $\Ind(G,g)\in \Gamma_2$. Since $P/P_1$ is a cyclic group and $P$ is an abelian group, we can write $P=\langle z, P_1\rangle$.

 We have $\Ind(G,g)\in \Gamma_1$ for any non-central $\{p,t\}'$-element $g$. Therefore for any $h$ such that $\Ind(G,h)\in \Gamma_2$ it is true that $C_G(h)/Z(G)$ is a $\{p,t\}$-group. In particular, $\Ind(G,h)_{\{p,t\}'}=|G||_{\{p,t\}'}$. Assume that there exists $z'\in P\setminus P_1$ such that $\Ind(G,z')\in \Gamma_1$. If $C_G(z')$ does not contain non-central $\{p,t\}'$-elements, then $\Ind(G,z)$ is connected to $\Ind(G,z')$ in $\Gamma( \Omega\setminus \{1\})$ and hence $\Ind(G,z')\in \Gamma_2$, contradicting $\Ind(G,z')\in \Gamma_1$. Let $s\in C_G(z')$ be a $\{p,t\}'$-element and $E\in Syl_p(C_G(z'))$. Since $C_G(x)$ contains some Sylow $p$-subgroup of $C_G(s)$, we can assume that $C_G(x^g)$ contains $E$ for some $g$. But $C_G(x^g)\cap P=P_1$, and we have a contradiction.

Let $y\in C_G(z)\setminus (Z(G)\cup P)$. As noted above, $y$ is a $t$-element. We can assume that $C_G(y)\cap P\in Syl_p(C_G(y))$. Hence $C_G(xy)\cap P\in Syl_p(C_G(xy))$. Obviously, $z$ and $P_1$ do not lie in $C_G(xy)$. Let $zg\in C_G(xy)$, where $g\in P_1\setminus Z(G)$. Let $\widetilde{\ }: G\rightarrow G/O$ be a natural homomorphism. Note $\widetilde{(xy)}=\widetilde{y}$. Hence $\widetilde{z}\in C_{\widetilde{G}}(\widetilde{xy})$, and thus $\widetilde{g}\in C_{\widetilde{G}}(\widetilde{xy })$. Since $|O|$ is coprime to $|g|$, then $\widetilde{C_G(g)}=C_{\widetilde{G}}(\widetilde{g})$. Hence $C_G(g)$ contains the group $O.\langle \widetilde{y}\rangle$, and therefore $y\in C_G(g)$ contradicting the fact that $C_G(y)\cap P<Z(G)$. Thus, it is proved that $P_1<Z(G)$.

Since $z$ acts without fixed points on $x^G$, $\Ind(G,z)>1$. Denote by $\Gamma'$ the connected component of $\Gamma(\Omega\setminus\{1\})$ containing $\Ind(G,z)$. Note that $\Ind(G,g)_{p'}\in\Gamma'\cup\{1\}$ for any $g\in C_G(z)$. Assume that there exists $h\not\in \langle z\rangle$ such that $\Ind(G,h)\in \Omega \setminus (\Gamma'\cup\{1\})$. Then $h$ centralizes some Sylow $p$-subgroup and, therefore, we can assume that $h\in C_G(z)$. Thus $\Ind(G,h)\in \Gamma'$, contradicting the hypotesis on $h$. We have $|\Gamma'|=1$ and $C_G(z)/Z(G)=\langle z\rangle$. Therefore $\Omega=N(\langle z\rangle)$, and in particular $\Gamma(\Omega)$ is connected, which is a contradiction.
\end{proof}

It follows from the proposition \ref{um} and \ref{dois} that $G\in R(p)^*$.

\begin{proposition}\label{tres}
If $G\in R(p)^*$ for some $p\in \pi(n)$ then $G=A\times B$, where $N(A)=\Omega$ and $N(B )=\{1,p^{\alpha}\}$. In particular, $n$ is a $p$-number.
\end{proposition}
\begin{proof}
Lemma \ref{star} implies that $G=N\leftthreetimes P$ where $P$ is a Sylow $p$-subgroup of $G$. Lemma \ref{starcor} implies that $Z(P)\leq Z(G)$.

Assume that there is $z\in P$ such that $\Ind(G,z)_{p'}>1$. The  separation of $\Omega$ implies that there exists $k\in\mu(\Omega)$ such that $k$ is not divisible by $\Ind(G,z)_{p'}$. Let $g\in G$ be such that $\Ind(G,g)=k$. We have $g=g_1g_2$, where $g_1$ is a $p'$-element and $g_2$ is a $p$-element. Since $C_G(g)=C_G(g_1)\cap C_G(g_2)$ and $\Ind(G,g)_p=1$, it follows that $\Ind(G,g_2)_p=1$. Hence $g_2\in Z(G)$. Thus, $\Ind(G,g)=\Ind(G,g_1)$. We have that $C_G(g_1)$ contains some Sylow $p$-subgroup of $G$, and therefore there is $z'\in C_G(g_1)\cap z^G$, deriving a contradiction.

Thus any $p$-element centralizes $N$ and hence $G\simeq N\times P$. Therefore for each $g\in P$ we have $\pi(\Ind(G,g))=\{p\}$. In particular $n$ is a $\{p\}$-number.
\end{proof}

The assertion of the theorem follows from Propositions \ref{um}, \ref{dois} and \ref{tres}.

\subsection*{Acknowledgment}
The work was supported by the grant of the President of the Russian Federation for
young scientists (MD-1264.2022.1.1).

\EditInfo{June 22, 2022}{January 11, 2023}{Ivan Kaygorodov}

\end{paper}